\documentclass[preprint,12pt]{elsarticle}
\usepackage{amssymb}
\usepackage{amsmath}

\newtheorem{thm}{Theorem}[section]
\newtheorem{lem}[thm]{Lemma}

\newdefinition{defn}[thm]{Definition}
\newdefinition{rem}[thm]{Remark}
\newproof{pf}{Proof}
\newproof{pot1}{Proof of Theorem \ref{jtj}}
\journal{arXiv}
\begin{document}
\begin{frontmatter}

\title{Ground State Solutions of Kirchhoff-type Fractional Dirichlet Problem with $p$-Laplacian}

\author{Taiyong Chen,\ \ Wenbin Liu
\footnote{Corresponding author.\\
{\it Telephone number:} (86-516) 83591530. {\it Fax number:} (86-516) 83591591.
{\it E-mail addresses:} taiyongchen@cumt.edu.cn (T. Chen), wblium@163.com (W. Liu), jinhua197927@163.com (H. Jin).}
,\ \ Hua Jin}

\address{Department of Mathematics, China University of Mining and Technology, Xuzhou 221116, PR China}

\begin{abstract}
We consider the Kirchhoff-type $p$-Laplacian Dirichlet problem containing the left and right fractional derivative operators. By using the Nehari method in critical point theory, we obtain the existence theorem of ground state solutions for such Dirichlet problem.
\end{abstract}


\end{frontmatter}

\section{Introduction}
\label{sec1}
In the present paper, we discuss the existence of ground state solutions for the Kirchhoff-type fractional Dirichlet problem with $p$-Laplacian of the form
\begin{eqnarray}
\label{kbvp}
\left\{
\begin{array}{ll}
\left(a+b\int_0^T|{_0D_t^\alpha}u(t)|^pdt\right)^{p-1}{_tD_T^\alpha}\phi_p({_0D_t^\alpha}u(t))=f(t,u(t)),\ \ t\in(0,T),\\
u(0)=u(T)=0,
\end{array}
\right.
\end{eqnarray}
where $a,b>0,\ p>1$ are constants, $_0D_t^\alpha$ and $_tD_T^\alpha$ are the left and right Riemann-Liouville fractional derivatives of order $\alpha\in(1/p,1]$ respectively, $\phi_p:\mathbb{R}\rightarrow\mathbb{R}$ is the $p$-Laplacian defined by
\begin{eqnarray*}
\phi_p(s)=|s|^{p-2}s\ (s\neq0),\ \ \phi_p(0)=0,
\end{eqnarray*}
and $f\in C^1([0,T]\times\mathbb{R},\mathbb{R})$.

The Kirchhoff equation (\cite{gy18}) is an extension of the wave equation which comes from the free vibrations of elastic strings, and takes into account the changes in length of the string produced by transverse vibrations. In addition, the fractional order models are more appropriate than the integer order models in real world owing to the fractional derivatives offer an wonderful tool to describe the memory and hereditary properties of a great deal of processes and materials (\cite{dkfa,3,hilf,kjfx,main}). Moreover the $p$-Laplacian (\cite{lsl}) often appears in non-Newtonian fluid theory, nonlinear elastic mechanics and so on.

Notice that, when $a=1$, $b=0$ and $p=2$, the left hand side of equation of BVP (\ref{kbvp}), which is nonlinear and nonlocal, reduces to the linear operator $_tD_T^\alpha{_0D_t^\alpha}$, and further reduces to the local operator $-d^2/dt^2$ when $\alpha=1$.

In recent years, there are many authors to study the fractional boundary value problems (BVPs for short) (\cite{21,zbhl,23,e1,25,wj}) and the Kirchhoff equations (\cite{gy2,gy4,gy6,gy11,gy19,gy25}), and obtain numerous important results. In addition, the models containing left and right fractional derivatives are recently gaining more attention (\cite{jbe3,jcre,jerw,12,fjy,zzr}) because of the applications in physical phenomena exhibiting anomalous diffusion.

Motivated by the above works, in this paper, we dicuss the existence of nontrivial ground state solutions for BVP (\ref{kbvp}). The main tool used here is the Nehari method.

For the nonlinearity $f$, we make the following assumptions throughout this paper.

$(H_{1})$ The mapping $x\rightarrow f(t,x)/|x|^{p^2-1}$ is strictly increasing on $\mathbb{R}\setminus\{0\}$ for $\forall t\in[0,T]$.

$(H_{2})$ $f(t,x)=o(|x|^{p-1})$ as $|x|\rightarrow0$ uniformly for $\forall t\in[0,T]$.

$(H_{3})$ There exist two constants $\mu>p^2,\ R>0$ such that
\begin{eqnarray*}
0<\mu F(t,x)\leq xf(t,x),\ \ \forall t\in[0,T],\ x\in\mathbb{R}\ \mbox{with}\ |x|\geq R,
\end{eqnarray*}
where $F(t,x)=\int_0^xf(t,s)ds$.

Now we state our main result.

\begin{thm}
\label{jtj}
Let $(H_{1})$-$(H_{3})$ be satisfied. Then BVP (\ref{kbvp}) possesses at least one nontrivial ground state solution.
\end{thm}

The rest of this paper is organized as follows. Some preliminary results are presented in Section \ref{sec2}. Section \ref{sec3} is devoted to prove Theorem \ref{jtj}.

\section{Preliminaries}
\label{sec2}

In this section, we present some basic definitions and notations of the fractional calculus (\cite{15,18}). Moreover we introduce a fractional Sobolev space and some properties of this space (\cite{fjy}).

\begin{defn}%[\cite{15}]
\label{defn2.1}
For $\gamma>0$, the left and right Riemann-Liouville fractional integrals of order $\gamma$ of a function $u:[a,b]\rightarrow\mathbb{R}$ are given by
\begin{eqnarray*}
&&_aI_{t}^\gamma u(t)=\frac{1}{\Gamma(\gamma)}\int_a^t(t-s)^{\gamma -1}u(s)ds,\\
&&_tI_{b}^\gamma u(t)=\frac{1}{\Gamma(\gamma)}\int_t^b(s-t)^{\gamma -1}u(s)ds,
\end{eqnarray*}
provided that the right-hand side integrals are pointwise defined on $[a,b]$, where $\Gamma(\cdot)$ is the Gamma function.
\end{defn}

\begin{defn}%[\cite{15}]
\label{defn2.2}
For $n-1\leq\gamma<n\ (n\in\mathbb{N})$, the left and right Riemann-Liouville fractional derivatives of order $\gamma$ of a function $u:[a,b]\rightarrow\mathbb{R}$ are given by
\begin{eqnarray*}
&&_aD_{t}^\gamma u(t)=\frac{d^n}{dt^n}{_a}I_{t}^{n-\gamma} u(t),\\
&&_tD_{b}^\gamma u(t)=(-1)^n\frac{d^n}{dt^n}{_t}I_{b}^{n-\gamma} u(t).
\end{eqnarray*}
\end{defn}

\begin{rem}
\label{zxzj}
When $\gamma=1$, one can obtain from Definition \ref{defn2.1} and \ref{defn2.2} that
\begin{eqnarray*}
_aD_{t}^1 u(t)=u'(t),\ \ _tD_{b}^1 u(t)=-u'(t),
\end{eqnarray*}
where $u'$ is the usual first-order derivative of $u$.
\end{rem}

\begin{defn}%[\cite{cjh}]
\label{defn3.1}
For $0<\alpha\leq1$ and $1<p<\infty$, the fractional derivative space $E{_0^{\alpha,p}}$ is defined by the closure of $C_0^\infty([0,T],\mathbb{R})$ with respect to the following norm
\begin{eqnarray*}
\|u\|_{E^{\alpha,p}}=(\|u\|_{L^p}^p+\|{_0D_t^\alpha}u\|_{L^p}^p)^{\frac{1}{p}},
\end{eqnarray*}
where $\|u\|_{L^p}=\left(\int_0^T|u(t)|^pdt\right)^{1/p}$ is the norm of $L^p([0,T],\mathbb{R})$.
\end{defn}

\begin{rem}
It is obvious that, for $u\in E{_0^{\alpha,p}}$, one has
\begin{eqnarray*}
u,{_0D_t^\alpha}u\in L^p([0,T],\mathbb{R}),\ \ u(0)=u(T)=0.
\end{eqnarray*}
\end{rem}

\begin{lem}[see \cite{fjy}]
\label{lem1}
Let $0<\alpha\leq1$ and $1<p<\infty$. The fractional derivative space $E{_0^{\alpha,p}}$ is a reflexive and separable Banach space.
\end{lem}

\begin{lem}[see \cite{fjy}]
Let $0<\alpha\leq1$ and $1<p<\infty$. For $u\in E_0^{\alpha,p}$, one has
\begin{eqnarray}
\label{cp}
\|u\|_{L^p}\leq C_p\|{_0D_t^\alpha}u\|_{L^p},
\end{eqnarray}
where
\begin{eqnarray*}
C_p=\frac{T^\alpha}{\Gamma(\alpha+1)}>0
\end{eqnarray*}
is a constant. Moreover, if $\alpha>1/p$, then
\begin{eqnarray}
\label{cwq}
\|u\|_\infty\leq C_\infty\|{_0D_t^\alpha}u\|_{L^p},
\end{eqnarray}
where $\|u\|_\infty=\max_{t\in[0,T]}|u(t)|$ is the norm of $C([0,T],\mathbb{R})$ and
\begin{eqnarray*}
C_\infty=\frac{T^{\alpha-\frac{1}p}}{\Gamma(\alpha)(\alpha q-q+1)^{\frac{1}q}}>0,\ \
q=\frac{p}{p-1}>1
\end{eqnarray*}
are two constants.
\end{lem}

\begin{rem}
By (\ref{cp}), we can consider the space $E_0^{\alpha,p}$ with norm
\begin{eqnarray}
\label{fsdj}
\|u\|_{E^{\alpha,p}}=\|{_0D_t^\alpha}u\|_{L^p}
\end{eqnarray}
in what follows.
\end{rem}

\begin{lem}[see \cite{fjy}]
\label{thm3.2}
Let $1/p<\alpha\leq1$ and $1<p<\infty$. The imbedding of $E_0^{\alpha,p}$ in $C([0,T],\mathbb{R})$ is compact.
\end{lem}

\section{Ground state solutions of BVP (\ref{kbvp})}
\label{sec3}

The purpose of this section is to prove our main result via the Nehari method. To this end, we are going to set up the corresponding variational framework of BVP (\ref{kbvp}).

Define the functional $I:E_0^{\alpha,p}\rightarrow\mathbb{R}$ by
\begin{align*}
%\label{fh}
I(u)
&=\frac{1}{bp^2}\left(a+b\int_0^T|{_0D_t^\alpha}u(t)|^pdt\right)^p-\int_0^TF(t,u(t))dt
-\frac{a^p}{bp^2}\nonumber\\
&=\frac{1}{bp^2}(a+b\|u\|_{E^{\alpha,p}}^p)^p-\int_0^TF(t,u(t))dt-\frac{a^p}{bp^2}.
\end{align*}
Then there is one-to-one correspondence between the critical points of energy functional $I$ and the weak solutions of BVP (\ref{kbvp}). It is easy to check from (\ref{cwq}), (\ref{fsdj}) and $f\in C^1([0,T]\times\mathbb{R},\mathbb{R})$ that the functional $I$ is well defined on $E_0^{\alpha,p}$ and is second-order continuously Fr\'{e}chet differentiable, that is, $I\in C^2(E_0^{\alpha,p},\mathbb{R})$. Furthermore we have
\begin{align*}
%\label{fhds}
\langle I'(u),v\rangle
&=(a+b\|u\|_{E^{\alpha,p}}^p)^{p-1}
\int_0^T\phi_p({_0D_t^\alpha}u(t)){_0D_t^\alpha}v(t)dt\nonumber\\
&\ \ \ \ -\int_0^Tf(t,u(t))v(t)dt,\ \ \forall u,v\in E_0^{\alpha,p},
\end{align*}
which yields
\begin{eqnarray*}
\label{fhds2}
\langle I'(u),u\rangle
=(a+b\|u\|_{E^{\alpha,p}}^p)^{p-1}\|u\|_{E^{\alpha,p}}^p-\int_0^Tf(t,u(t))u(t)dt.
\end{eqnarray*}

Now let us define
\begin{eqnarray*}
\label{nlx}
\mathcal{N}=\{u\in E_0^{\alpha,p}\setminus\{0\}|G(u)=0\},
\end{eqnarray*}
where
\begin{eqnarray*}
G(u)=\langle I'(u),u\rangle.
\end{eqnarray*}
Thus we know that any non-zero critical point of $I$ must be on $\mathcal{N}$. In the following, for simplicity, let
\begin{eqnarray*}
M_u=a+b\|u\|_{E^{\alpha,p}}^p.
\end{eqnarray*}
From $(H_1)$, one has
\begin{eqnarray}
\label{h1}
f'_2(t,x)x^2>(p^2-1)f(t,x)x,\ \ \forall (t,x)\in[0,T]\times(\mathbb{R}\setminus\{0\}),
\end{eqnarray}
where $f'_2(t,x)=\frac{\partial f(t,x)}{\partial x}$. Then, for $u\in\mathcal{N}$, we have
\begin{align}
\label{c1lx}
\langle G'(u),u\rangle
&=bp(p-1)M_u^{p-2}\|u\|_{E^{\alpha,p}}^{2p}+pM_u^{p-1}\|u\|_{E^{\alpha,p}}^p\nonumber\\
&\ \ \ \ -\int_{0}^{T}f'_2(t,u(t))u^2(t)dt-\int_{0}^{T}f(t,u(t))u(t)dt\nonumber\\
&<M_u^{p-2}\|u\|_{E^{\alpha,p}}^p(bp^2\|u\|_{E^{\alpha,p}}^p+ap)-p^2\int_{0}^{T}f(t,u(t))u(t)dt\nonumber\\
&=a(p-p^2)M_u^{p-2}\|u\|_{E^{\alpha,p}}^p\leq0,
\end{align}
which means that $\mathcal{N}$ has a $C^1$ structure and is a manifold.

\begin{lem}
\label{zrys}
Assume $(H_{1})$ holds. If $u\in\mathcal{N}$ is a critical point of $I|_\mathcal{N}$, then $I'(u)=0$, that is, $\mathcal{N}$ is a natural constraint for $I$.
\end{lem}

\begin{pf}
If $u\in\mathcal{N}$ is a critical point of $I|_\mathcal{N}$, then there exists a Lagrange multiplier $\lambda\in\mathbb{R}$ such that
\begin{eqnarray*}
I'(u)=\lambda G'(u).
\end{eqnarray*}
Then we get
\begin{eqnarray*}
\langle I'(u),u\rangle=\lambda\langle G'(u),u\rangle=0,
\end{eqnarray*}
which together with (\ref{c1lx}) yields $\lambda=0$. So we have $I'(u)=0$. $\Box$
\end{pf}

In order to discuss the critical points of $I|_\mathcal{N}$, we need to investigate the structure of $\mathcal{N}$.

\begin{lem}
\label{czwy}
Assume $(H_{1})$-$(H_{3})$ hold. For each $u\in E_0^{\alpha,p}\setminus\{0\}$, there is a unique $s=s(u)\in\mathbb{R}^+$ such that $su\in\mathcal{N}$.
\end{lem}

\begin{pf}
First, we claim that there exist constants $\rho,\sigma>0$ such that
\begin{eqnarray}
\label{fpf}
I(u)>0, \forall u\in B_\rho(0)\setminus\{0\},\ \ I(u)\geq\sigma, \forall u\in\partial B_\rho(0),
\end{eqnarray}
where $B_\rho(0)$ is an open ball in $E_0^{\alpha,p}$ with the radius $\rho$ and centered at $0$, and $\partial B_\rho(0)$ denote its boundary. That is, by $I(0)=0$, $0$ is a strict local minimizer of $I$. In fact, from $(H_{2})$, there are two constants $0<\varepsilon<1,\delta>0$ such that
\begin{eqnarray}
\label{530.1}
F(t,x)\leq\frac{(1-\varepsilon)a^{p-1}}{pC_p^p}|x|^p,
\ \ \forall (t,x)\in[0,T]\times[-\delta,\delta],
\end{eqnarray}
where $C_p>0$ is a constant defined in (\ref{cp}). Let $\rho=\delta/C_\infty$ and $\sigma=\varepsilon a^{p-1}\rho^p/p$, where $C_\infty>0$ is a constant defined in (\ref{cwq}). Then, by (\ref{cwq}) and (\ref{fsdj}), one has
\begin{eqnarray*}
\|u\|_\infty\leq C_\infty\|u\|_{E^{\alpha,p}}\leq\delta,\ \
\forall u\in \overline{B_\rho(0)},
\end{eqnarray*}
which together with (\ref{cp}), (\ref{fsdj}) and (\ref{530.1}) yields
\begin{align*}
I(u)&=\frac{1}{bp^2}M_u^p-\int_0^TF(t,u(t))dt-\frac{a^p}{bp^2}\\
&\geq\frac{a^{p-1}}{p}\|u\|_{E^{\alpha,p}}^p
-\frac{(1-\varepsilon)a^{p-1}}{pC_p^p}\int_0^T|u(t)|^pdt\\
&\geq\frac{a^{p-1}}{p}\|u\|_{E^{\alpha,p}}^p
-\frac{(1-\varepsilon)a^{p-1}}{p}\|u\|_{E^{\alpha,p}}^p\\
&=\frac{\varepsilon a^{p-1}}{p}\|u\|_{E^{\alpha,p}}^p=\sigma,\ \ \forall u\in \partial B_\rho(0).
\end{align*}

Second, we claim that $I(\xi u)\rightarrow-\infty$ as $\xi\rightarrow\infty$. In fact, from $(H_{3})$, a simple argument can show that there are two constants $c_1,c_2>0$ such that
\begin{eqnarray*}
F(t,x)\geq c_1|x|^\mu-c_2,\ \ \forall (t,x)\in[0,T]\times\mathbb{R}.
\end{eqnarray*}
Thus, for each $u\in E_0^{\alpha,p}\setminus\{0\},\ \xi\in\mathbb{R}^+$, we obtain from $\mu>p^2$ that
\begin{align*}
I(\xi u)&=\frac{1}{bp^2}M_{\xi u}^p-\int_0^TF(t,\xi u(t))dt-\frac{a^p}{bp^2}\\
&\leq\frac{1}{bp^2}M_{\xi u}^p-c_1\int_0^T|\xi u(t)|^\mu dt+c_2T-\frac{a^p}{bp^2}\\
&=\frac{1}{bp^2}(a+b\xi^p\|u\|_{E^{\alpha,p}}^p)^p
-c_1\xi^\mu\|u\|_{L^\mu}^\mu+c_2T-\frac{a^p}{bp^2}\\
&\rightarrow-\infty\ \ \mbox{as}\ \xi\rightarrow\infty.
\end{align*}

Let
\begin{eqnarray*}
g_u(s)=I(su),\ \ \forall s\in\mathbb{R}^+.
\end{eqnarray*}
Then, from what we have proved, $g_u$ has at least one maximum point $s(u)$ with maximum value greater than $\sigma>0$. Next, we prove that $g_u$ has a unique critical point for $s\in\mathbb{R}^+$, which then must be the global maximum point. Consider a critical point of $g_u$, one has
\begin{align*}
g_u'(s)
&=\langle I'(su),u\rangle\\
&=\|u\|_{E^{\alpha,p}}^pM_{su}^{p-1}s^{p-1}-\int_{0}^{T}f(t,su(t))u(t)dt\\
&=0,
\end{align*}
which together with (\ref{h1}) yields
\begin{align}
\label{wyx}
g_u''(s)
&=bp(p-1)\|u\|_{E^{\alpha,p}}^{2p}M_{su}^{p-2}s^{2p-2}\nonumber\\
&\ \ \ \ +(p-1)\|u\|_{E^{\alpha,p}}^pM_{su}^{p-1}s^{p-2}-\int_{0}^{T}f_2'(t,su(t))u^2(t)dt\nonumber\\
&<\|u\|_{E^{\alpha,p}}^pM_{su}^{p-2}\left(b(p^2-1)\|u\|_{E^{\alpha,p}}^ps^{2p-2}
+a(p-1)s^{p-2}\right)\nonumber\\
&\ \ \ \ -\frac{p^2-1}{s}\int_{0}^{T}f(t,su(t))u(t)dt\nonumber\\
&=\|u\|_{E^{\alpha,p}}^pM_{su}^{p-2}\left(b(p^2-1)\|u\|_{E^{\alpha,p}}^ps^{2p-2}
+a(p-1)s^{p-2}\right)\nonumber\\
&\ \ \ \ -(p^2-1)\|u\|_{E^{\alpha,p}}^pM_{su}^{p-1}s^{p-2}\nonumber\\
&=a\|u\|_{E^{\alpha,p}}^pM_{su}^{p-2}(p-p^2)s^{p-2}\leq0.
\end{align}
Hence, if $s$ is a critical point of $g_u$, then it must be a strict local maximum point. This ensures the uniqueness of critical point of $g_u$. Finally, from
\begin{eqnarray}
\label{guds}
g_u'(s)=\frac{1}{s}\langle I'(su),su\rangle,\ \ \forall t\in\mathbb{R}^+,
\end{eqnarray}
we obtain that, if $s$ is a critical point of $g_u$, then $su\in\mathcal{N}$. $\Box$
\end{pf}

Let us define
\begin{eqnarray*}
m=\inf_\mathcal{N}I.
\end{eqnarray*}
Then we get from (\ref{fpf}) that
\begin{eqnarray*}
m\geq\inf_{\partial B_\rho(0)}I\geq\sigma>0.
\end{eqnarray*}

\begin{lem}
\label{kd}
Assume $(H_{1})$-$(H_{3})$ hold. Then there exists $u^*\in\mathcal{N}$ such that $I(u^*)=m$.
\end{lem}

\begin{pf}
By Lemma \ref{thm3.2}, we obtain that the functional
\begin{eqnarray*}
u\rightarrow\int_0^TF(t,u(t))dt,\ \ \forall u\in E_0^{\alpha,p}
\end{eqnarray*}
is weakly continuous. Thus, as the sum of a convex continuous functional and a weakly continuous one, $I$ is weakly lower semi-continuous on $E_0^{\alpha,p}$.

Let $\{u_k\}\subset\mathcal{N}$ be a minimizing sequence of $I$, then one has
\begin{eqnarray}
\label{wqxl}
I(u_k)=m+o(1),\ \ G(u_k)=0.
\end{eqnarray}
Next, we prove that $\{u_k\}$ is bounded in $E_0^{\alpha,p}$. Based on the continuity of $\mu F(t,x)-xf(t,x)$ and $(H_{3})$, we see that there exists a constant $c>0$ such that
\begin{eqnarray*}
F(t,x)\leq\frac{1}{\mu}xf(t,x)+c,\ \ \forall (t,x)\in[0,T]\times\mathbb{R}.
\end{eqnarray*}
Thus, from (\ref{wqxl}), we have
\begin{align*}
m+o(1)&=I(u_k)\\
&\geq\frac{1}{bp^2}M_{u_k}^p
-\frac{1}{\mu}\int_0^Tf(t,u_k(t))u_k(t)dt-cT-\frac{a^p}{bp^2}\\
&=\frac{1}{bp^2}M_{u_k}^p
-\frac{1}{\mu}M_{u_k}^{p-1}\|u_k\|_{E^{\alpha,p}}^p-cT-\frac{a^p}{bp^2}\\
&=M_{u_k}^{p-1}\left(\left(\frac{1}{p^2}-\frac{1}{\mu}\right)
\|u_k\|_{E^{\alpha,p}}^p+\frac{a}{bp^2}\right)-cT-\frac{a^p}{bp^2}.
\end{align*}
Hence it follows from $\mu>p^2$ that $\{u_k\}$ is bounded in $E_0^{\alpha,p}$.

Since $E_0^{\alpha,p}$ is a reflexive Banach space (see Lemma \ref{lem1}), up to a subsequence, we can assume $u_k\rightharpoonup u$ in $E_0^{\alpha,p}$. Moreover, from Lemma \ref{thm3.2}, one has $u_k\rightarrow u$ in $C([0,T],\mathbb{R})$. Next, we prove $u\neq0$. By $(H_2)$, we get that for $\forall\varepsilon>0$, there exists a constant $\delta>0$ such that
\begin{eqnarray*}
f(t,x)x\leq\varepsilon|x|^p,\ \ \forall (t,x)\in[0,T]\times[-\delta,\delta].
\end{eqnarray*}
Then, assume $\|u_k\|_\infty\leq\delta$, we obtain from (\ref{cwq}), (\ref{fsdj}) and $u_k\in\mathcal{N}$ that
\begin{align*}
C_\infty^{-p}(a+bC_\infty^{-p}\|u_k\|_\infty^p)^{p-1}\|u_k\|_\infty^p
&\leq(a+b\|u_k\|_{E^{\alpha,p}}^p)^{p-1}\|u_k\|_{E^{\alpha,p}}^p\\
&=\int_0^Tf(t,u_k(t))u_k(t)dt\\
&\leq\varepsilon\int_0^T|u_k(t)|^pdt\\
&\leq\varepsilon T\|u_k\|_\infty^p,
\end{align*}
which is a contradiction. Hence we have
\begin{eqnarray*}
\|u\|_\infty=\lim_{k\rightarrow\infty}\|u_k\|_\infty\geq\delta>0,
\end{eqnarray*}
and then $u\neq0$. Thus, by Lemma \ref{czwy}, there exists $s\in\mathbb{R}^+$ such that $su\in\mathcal{N}$. Therefore, together with the fact that $I$ is weakly lower semi-continuous, we obtain
\begin{eqnarray}
\label{jxbds}
m\leq I(su)\leq\varliminf_{k\rightarrow\infty}I(su_k)\leq\lim_{k\rightarrow\infty}I(su_k).
\end{eqnarray}

Finally, for $\forall u_k\in\mathcal{N}$, we see from (\ref{wyx}) and (\ref{guds}) that $s=1$ is the global maximum point of $g_{u_k}$. So one has
\begin{eqnarray*}
I(su_k)\leq I(u_k),
\end{eqnarray*}
which together with (\ref{jxbds}) implies
\begin{eqnarray*}
m\leq I(su)\leq\lim_{k\rightarrow\infty}I(u_k)=m.
\end{eqnarray*}
That is, $m$ is achieved at $su\in\mathcal{N}$. $\Box$
\end{pf}

Now we give the proof of our main result.

\begin{pot1}
By Lemma \ref{kd}, we get $u^*\in\mathcal{N}$ such that $I(u^*)=m=\inf_\mathcal{N}I>0$, that is, $u^*$ is a non-zero critical point of $I|_\mathcal{N}$. Then, from Lemma \ref{zrys}, we know $I'(u^*)=0$, and so $u^*$ is a nontrivial ground state solution of BVP (\ref{kbvp}). $\Box$
\end{pot1}

%\section*{Competing interests}
%The authors declare that they have no competing interests.

%\section*{Authors' contributions}
%The authors contributed equally in this article. They read and approved the final manuscript.

\section*{Acknowledgements}
This work was supported by the National Natural Science Foundation of China (11271364) and the Nature Science Foundation of Jiangsu Province (BK20130170).


\section*{References}
\begin{thebibliography}{0000}

\bibitem{21} R.P. Agarwal, D. O'Regan, S. Stanek,
Positive solutions for Dirichlet problems of singular nonlinear fractional differential equations,
{\it J. Math. Anal. Appl.} 371 (2010) 57-68.

\bibitem{gy2} A. Arosio, S. Panizzi,
On the well-posedness of the Kirchhoff string,
{\it Trans. Amer. Math. Soc.} 348 (1996) 305-330.

\bibitem{zbhl} Z. Bai, H. L\"{u},
Positive solutions for boundary value problem of nonlinear fractional differential equation,
{\it J. Math. Anal. Appl.} 311 (2005) 495-505.

\bibitem{23} M. Benchohra, S. Hamani, S.K. Ntouyas,
Boundary value problems for differential equations with fractional order and nonlocal conditions,
{\it Nonlinear Anal.} 71 (2009) 2391-2396.

\bibitem{jbe3} D.A. Benson, S.W. Wheatcraft, M.M. Meerschaert,
The fractional-order governing equation of L\'{e}vy motion,
{\it Water Resour. Res.} 36 (2000) 1413-1423.

\bibitem{gy4} S. Bernstein,
Sur une classe d'\'{e}quations fonctionnelles aux d\'{e}iv\'{e}s partielles,
{\it Bull. Acad. Sci. URSS. S\'{e}. Math.} 4 (1940) 17-26.

\bibitem{e1} G.M. Bisci, D. Repovs,
Higher nonlocal problems with bounded potential,
{\it J. Math. Anal. Appl.} 420 (2014) 167-176.

\bibitem{gy6} M.M. Cavalcanti, V.N. Domingos Cavalcanti, J.A. Soriano,
Global existence and uniform decay rates for the Kirchhoff-Carrier equation with nonlinear dissipation,
{\it Adv. Differential Equations} 6 (2001) 701-730.

\bibitem{jcre} J. Cresson,
Inverse problem of fractional calculus of variations for partial differential equations,
{\it Commun. Nonlin. Sci. Numer. Simul.} 15 (2010) 987-996.

\bibitem{gy11} P. D'Ancona, S. Spagnolo,
Global solvability for the degenerate Kirchhoff equation with real analytic data,
{\it Invent. Math.} 108 (1992) 247-262.

\bibitem{25} M.A. Darwish, S.K. Ntouyas,
On initial and boundary value problems for fractional order mixed type functional differential inclusions,
{\it Comput. Math. Appl.} 59 (2010) 1253-1265.

\bibitem{dkfa} K. Diethelm, A.D. Freed,
On the solution of nonlinear fractional order differential equations used in the modeling of viscoelasticity, in: F. Keil, W. Mackens, H. Voss, J. Werther (Eds.), Scientific Computing in Chemical Engineering II-Computational Fluid Dynamics, Reaction Engineering and Molecular Properties,
Springer-Verlag, Heidelberg, 1999, 217-224.

\bibitem{jerw} V.J. Ervin, J.P. Roop,
Variational formulation for the stationary fractional advection dispersion equation,
{\it Numer. Meth. Part. Diff. Eqs.} 22 (2006) 558-576.

\bibitem{12} G.J. Fix, J.P. Roop,
Least squares finite-element solution of a fractional order two-point boundary value problem,
{\it Comput. Math. Appl.} 48 (2004) 1017-1033.

\bibitem{3} W.G. Glockle, T.F. Nonnenmacher,
A fractional calulus approach of self-similar protein dynamcs,
{\it Biophys. J.} 68 (1995) 46-53.

\bibitem{hilf} R. Hilfer,
Applications of Fractional Calculus in Physics,
World Scientific, Singapore, 2000.

\bibitem{wj} W. Jiang,
The existence of solutions to boundary value problems of fractional differential equations at resonance,
{\it Nonlinear Anal.} 74 (2011) 1987-1994.

\bibitem{fjy} F. Jiao, Y. Zhou,
Existence results for fractional boundary value problem via critical point theory,
{\it Internat. J. Bifur. Chaos} 22 (2012) Article ID 1250086.

\bibitem{15} A.A. Kilbas, H.M. Srivastava, J.J. Trujillo,
Theory and Applications of Fractional Differential Equations,
Elsevier, Amsterdam, 2006.

\bibitem{gy18} G. Kirchhoff, Mechanik, Teubner, Leipzig, 1883.

\bibitem{kjfx} J.W. Kirchner, X. Feng, C. Neal,
Fractal stream chemistry and its implications for contaminant transport in catchments,
{\it Nature} 403 (2000) 524-526.

\bibitem{lsl} L.S. Leibenson,
General problem of the movement of a compressible fluid in a porous medium,
{\it Izvestiia Akademii Nauk Kirgizsko\u{\i} SSSR} 9 (1983) 7-10.

\bibitem{gy19} J.L. Lions,
On some questions in boundary value problems of mathematical physics, in: Contemporary Developments in Continuum Mechanics and Partial Differential Equations, Proc. Internat. Sympos., Inst. Mat., Univ. Fed. Rio de Janeiro, Rio de Janeiro, 1977, in: North-Holland Math. Stud., vol. 30, North-Holland, Amsterdam, New York, 1978, 284-346.

\bibitem{main} F. Mainardi,
Fractional calculus: some basic problems in continuum and statistical mechanics,
in: A. Carpinteri, F. Mainardi (Eds.), Fractals and Fractional Calculus in Continuum Mechanics,
Springer-Verlag, Wien, 1997, 291-348.

%\bibitem{dl1} J. Mawhin, M. Willem,
%Critical point theory and Hamiltonian systems,
%Applied Mathematical Sciences, Springer, New York, 1989.

\bibitem{gy25} S.I. Pohoz\u{a}ev,
A certain class of quasilinear hyperbolic equations,
{\it Mat. Sb. (NS)} 96 (1975) 152-166.

%\bibitem{dl2} P. Rabinowitz,
%Minimax methods in critical point theory with applications to differential equations,
%in: CBMS Regional Conference Series in Mathematics, American Mathematical Society, Providence, RI, 1986.

\bibitem{18} S.G. Samko, A.A. kilbas, O.I. Marichev,
Fractional Integrals and Detivatives: Theory and Applications,
Gordon and Breach, New York, 1993.

%\bibitem{js} J. Simon,
%R\'{e}gularit\'{e} de la solution d'un probl\`{e}me aux limites non lin\'{e}aires,
%{\it Ann. Fac. Sci. Tolouse} 3 (1981) 247-274.

\bibitem{zzr} Z. Zhang, R. Yuan,
Infinitely-many solutions for subquadratic fractional Hamiltonian systems with potential changing sign,
{\it Adv. Nonlinear Anal.} 4 (2015) 59-72.

\end{thebibliography}
\end{document}